\documentclass[12pt,a4paper, twoside, reqno]{amsart} %Change as appropriate
\usepackage{graphicx}
\usepackage[T1]{fontenc} % Font encoding.
\usepackage[latin1]{inputenc} % Input encoding - change for other
\usepackage{bbm} % More fonts 
\usepackage[english]{babel} % Multiple languages support.
\usepackage{varioref} % A few fancy references 

\usepackage{amsmath,amssymb,mathrsfs,amsthm,amsthm}
\usepackage{setspace} 
\begin{document}

\title[Tensor Product Factorization]{A Tensor Product Factorization for Certain Tilting Modules}
\author[M. Fazeel Anwar]{M. Fazeel Anwar}
%(Received 07 January 2010)}}

\footnotetext[1]{Research supported by COMSATS Institute of Information Technology (CIIT), Islamabad, Pakistan.}

\begin{abstract}

\noindent Let $G$ be a semisimple, simply connected linear algebraic group over an algebraically closed field $k$ of characteristic $p>0$. In a recent paper \cite{Doty 2009}, Doty introduces the notion of $r-$minuscule weight and exhibits a tensor product factorization of a corresponding tilting module under the assumption $p \geq 2h-2$, where $h$ is the Coxeter number. We remove this restriction and consider some variations involving the more general notion of $(p,r)-$minuscule weights.

\bigskip


\bigskip


\end{abstract}

\maketitle

\noindent Let $G$ be a semisimple, simply connected linear algebraic group over an algebraically closed field $k$ of characteristic $p>0$. In a recent paper \cite{Doty 2009}, Doty observed that the tensor product of the Steinberg module with a minuscule module is always indecomposable tilting. In this paper, we show that the tensor product of the Steinberg module with a module whose dominant weights are $p-$minuscule is a tilting module, not always indecomposable. We also give the decomposition of such a module into indecomposable tilting modules. Doty also proved that if $p \geq 2h-2$, then for $r-$minuscule weights the tilting module is isomorphic to a tensor product of two simple modules, usually in many ways. We remove the characteristic restriction on this result. A generalization of \cite[proposition 5.5(i)]{Donkin 2007} for $(p,r)-$ minuscule weights is also given. We start by setting up notation and stating some important definitions and results which will be useful later on.

\bigskip

\noindent Let $F: G \rightarrow G$ be the Frobenius morphism of $G$. Let $B$ be a Borel subgroup of $G$ and $T \subset B$ be a maximal torus of $G$. Let ${\rm mod}(G)$ be the category of finite dimensional rational $G-$modules. Define $X(T)$ to be the group of multiplication characters of $T$. For a $T-$module $V$ and $\lambda \in X(T)$, write $V^{\lambda}$ for corresponding weight space of $V$. Those $\lambda$'s for which $V^{\lambda}$ is non-zero are called weights of $V$. Any $G-$module $M$ is completely reducible as a $T-$module. So $M$ decomposes as a direct sum of its weight spaces and we have $M=\oplus_{\lambda \in X(T)} \,M^{\lambda}$ as a $T-$module. We will write $M^{[1]}$ for $M^{F}$. The Weyl group $W$ acts on $T$ in the usual way. Let $ \{e(\lambda), \lambda \in X(T)\}$ be the canonical basis for the integral group ring $\mathbb{Z}X(T)$. The character of $M$ is defined by ${\rm ch}\,M=\Sigma_{\lambda \in X(T)}({\rm dim}\,M^{\lambda})\,e(\lambda)$. For $\phi = \Sigma_{\mu \in X(T)}\, a_{\mu} e(\mu) \in \mathbb{Z}X(T)$ we set $\phi^{[1]} = \Sigma_{\mu \in X(T)}\, a_{\mu} e(p\mu) \in \mathbb{Z}X(T)$. Let $\Phi$ be the system of roots, $\Phi^{+}$ the system of positive roots that make $B$ the negative Borel. Let $S$ be the set of simple roots. For $\alpha \in \Phi$ the co-root of $\alpha$ is denoted by $\alpha^{v}$. Let $X^{+}(T)$ denotes the set of dominant weights. For $\lambda \in X^{+}(T)$, we define the set of $r-$restricted weights $X_{r}(T)$ by $$ X_{r}(T)=\{ (\lambda,\alpha^{v})<p^{r} : \text{for all simple roots } \alpha \}.$$ 

\bigskip

\noindent For $\lambda \in X^{+}(T)$, let $\Delta(\lambda)$ denote the Weyl module of highest weight $\lambda$, the dual Weyl module of highest weight $\lambda$ is denoted by $\nabla(\lambda)$, and $L(\lambda)$ denotes the simple module of highest weight $\lambda$. The dual Weyl module $\nabla(\lambda)$ has simple socle $L(\lambda)$ and $\{L(\lambda) \text{ : }\lambda \in X^{+}(T)\}$ is a complete set of pairwise non-isomorphic simple $G-$modules. A good filtration of a $G-$module $M$ is defined as a filtration $0=M_{0} \leq M_{1} \leq M_{2} \leq  ... \leq M_{n} =M$ such that for each $0 <i \leq n$, $M_{i}/M_{i-1}$ is either zero or isomorphic to $\nabla(\lambda_{i})$ for some $\lambda_{i} \in X^{+}(T)$. 

\bigskip

\noindent A tilting module of $G$ is a finite dimensional $G-$module $M$ such that $M$ and its dual module $M^{*}$ both admit good filtrations. For each $ \lambda \in X^{+}(T)$ there is an indecomposable tilting module $T(\lambda)$ which has highest weight $\lambda$. Every tilting module is a direct sum of copies of $T(\lambda),\, \lambda \in X^{+}(T)$. For $\lambda \in X^{+}(T)$ the tilting module $T((p-1)\rho+\lambda)$ is projective as a $G_{1}-$module, where $G_{1}$ is the first infinitesimal group and $\rho$ is the half sum of positive roots.  

\bigskip

\noindent Now for $\lambda \in X_{r}(T)$, the modules $L(\lambda)$ form a complete set of pairwise non-isomorphic irreducible $G_{r}$ modules. For $\mu \in X(T)$ let $\hat{Q}_{r}(\mu)$ denote the projective cover of $L(\mu)$ as a $G_{r}T$ module see e.g \cite{Jantzen} and \cite{Jantzen 1979}. The modules $\hat{Q}_{r}(\lambda),\, \lambda \in X_{r}(T)$, form a complete set of pairwise non-isomorphic projective $G_{r}$ modules. We refer the reader to \cite{Jantzen}, \cite{Donkin 2007} and \cite{Donkin 1980} for terminology and results not explained here.

\bigskip

\noindent A dominant weight $\lambda$ is called minuscule if the weights of $\Delta(\lambda)$ form a single orbit under the action of $W$. Equivalently, by \cite[chapter VIII, Section 7, proposition 6(iii)]{Bourbaki 1975}, $\lambda$ is minuscule if $-1 \leq (\lambda, \alpha^{v}) \leq 1$ for all roots $\alpha$. If $s(\lambda)=\Sigma_{\mu \in W\lambda} \, e(\mu)$ then $\lambda$ minuscule implies $s(\lambda)={\rm ch}\,\Delta(\lambda)={\rm ch}\,\nabla(\lambda)={\rm ch}\,L(\lambda)$. For $\lambda \in X^{+}(T)$ define $\lambda$ to be $p-$minuscule if $\langle \lambda,\beta_{0}^{v} \rangle \leq p$, where $\beta_{0}$ is the highest short root. Moreover we define a weight $\lambda \in X_{r}(T)$ to be $(p,r)$-minuscule if $\lambda=\Sigma_{j=0}^{r-1}\, p^{j} \lambda^{j}$, where each $\lambda^{j}$ is $p-$minuscule (and $\lambda^{j} \in X_{1}(T)$). In \cite{Doty 2009} Doty defines a weight $\lambda$ to be $r-$minuscule if $\lambda=\Sigma_{j=0}^{r-1}\, p^{j} \lambda^{j}$, with each $\lambda^{j}$ minuscule. Note that $\lambda$ minuscule implies $\lambda$ is $p-$minuscule. Similarly if $\lambda$ is $r-$minuscule then $\lambda$ is $(p,r)-$minuscule.

\bigskip

\noindent {\bf Definition.} For $\lambda=\Sigma_{j=0}^{r-1}\, p^{j} \lambda^{j} \, \in X_{r}(T)$, $\lambda^{j} \in X_{1}(T)$ define $$s_{r}(\lambda)=s(\lambda^{0})s(p\lambda^{1}) \, ... \,s(p^r\lambda^{r-1}).$$

\bigskip

\noindent {\bf Proposition 1.} If $\lambda$ is $(p,r)-$minuscule then $${\rm ch}\,T((p^{r}-1)\rho+\lambda)=\chi((p^{r}-1)\rho)s_{r}(\lambda).$$

\bigskip

\noindent {\bf Proof.} By \cite[theorem 5.3]{Donkin 2007} we have if $\lambda \in X_{1}(T)$ and $T((p-1)\rho+\lambda)|_{G_{1}}$ is indecomposable then $T((p-1)\rho+\lambda) \otimes T(\mu)^{[1]} \simeq T((p-1)\rho+\lambda+p\mu)$ for all $\mu \in X^{+}(T)$. Also by the argument of \cite[proposition 5.5]{Donkin 2007} for $p-$minuscule (and restricted) $\lambda$ we get that $T((p-1)\rho+\lambda)|_{G_{1}}$ is indecomposable. So we have $T((p^{r}-1)\rho+\lambda)=\bigotimes_{j=0}^{r-1}\,T((p-1)\rho+\lambda^{j})^{[j]}$. So ${\rm ch}\,T((p^{r}-1)\rho+\lambda)=\Pi_{j=0}^{r-1}\,{\rm ch}\,T((p-1)\rho+\lambda^{j})^{[j]}$. Since each $\lambda^{j}$ is $p-$minuscule by \cite[proposition 5.5]{Donkin 2007} we get ${\rm ch}\,T((p-1)\rho+\lambda^{j})= \chi((p-1)\rho)s(\lambda^{j})$. Hence ${\rm ch}\,T((p^{r}-1)\rho+\lambda)=\Pi_{j=0}^{r-1}\, \chi((p-1)\rho)^{[j]}s(\lambda^{j})^{[j]}$. Combine this with above definition to get the result.

\bigskip

\noindent \textit{Remark.} If $\lambda$ is minuscule then $s(\lambda)={\rm ch}\, L(\lambda)$ and hence $T((p-1)\rho+\lambda)={\rm St} \otimes L(\lambda)$. Because these are tilting modules with same character. This gives us \cite[lemma]{Doty 2009}.

\bigskip

\bigskip

\noindent {\bf Lemma 1.}
 
\noindent {\bf (a)} if $\mu \in X^{+}(T)$ then $T((p^{r}-1)\rho)\otimes T(\mu)^{[r]} \simeq T((p^{r}-1)\rho+p^{r}\mu)$.

\bigskip

\noindent {\bf (b)} suppose $\lambda$ is minuscule then ${\rm St} \otimes L(\lambda) \simeq \hat{Q}_1((p-1)\rho+w_0\lambda)$ as $G_1T$ modules, where $w_{0}$ is the longest element of $W$. In particular ${\rm St} \otimes L(\lambda)|_{G_{1}}$ is indecomposable.

\bigskip

\noindent {\bf (c)} if $\lambda$ is minuscule and $\mu \in X^{+}(T)$ then $${\rm St} \otimes L(\lambda) \otimes T(\mu)^{[r]} \simeq T((p-1)\rho+\lambda+p^{r}\mu).$$

\bigskip

\noindent {\bf Proof.}  By \cite[II, 3.19]{Jantzen} with $i=0$ we have ${\rm St}_{r} \otimes \nabla(\mu)^{[r]} \simeq \nabla((p^{r}-1)\rho+p^{n}\mu)$ for every $\mu \in X^{+}(T)$. It follows that ${\rm St}_{r} \otimes V^{[r]}$ is tilting for every tilting module $V$. In particular ${\rm St}_{r} \otimes T(\mu)^{[r]}$ is tilting. By \cite[2.1]{Donkin 1993}, ${\rm St}_{r} \otimes T(\mu)^{[r]}$ is isomorphic to $T((p^{r}-1)\rho+p^{r}\mu)$. This proves part (a).

\bigskip

\noindent Since 
\begin{align*}
{\rm Hom}_{G_1T}(&L((p-1)\rho+w_{0}\lambda),{\rm St} \otimes L(\lambda))\\&= {\rm Hom}_{G_1T}(L((p-1)\rho+w_{0}\lambda)\otimes L(\lambda)^{*},{\rm St})\\&={\rm Hom}_{G_1T}(L((p-1)\rho+w_{0}\lambda)\otimes L(-w_{0}\lambda),{\rm St})\neq 0.
\end{align*}

\noindent we have $${\rm St} \otimes L(\lambda)|_{G_{1}}= \hat{Q}_1((p-1)\rho+w_0\lambda)\oplus Z.$$ Also by \cite[1.2(2)]{Donkin 2007JA}, ${\rm ch}\,\hat{Q}_1((p-1)\rho+w_0\lambda)=\chi((p-1)\rho)\psi$, where $\psi=\Sigma a_{\xi}e(\xi)$ and $a_{\xi} \geq 0$ for all $\xi$. 

\noindent Also by \cite[II, 11.7, lemma(c)]{Jantzen}, ${\rm ch}\,\hat{Q}_1((p-1)\rho+w_0\lambda)$ is $W$ invariant. This implies $\psi$ is $W$ invariant. Moreover $\hat{Q}_1((p-1)\rho+w_0\lambda)$ has unique highest weight $(p-1)\rho+\lambda$, so $\psi=s(\lambda)+\text{lower terms}$. But $\psi$ is $W$ invariant and ${\rm ch}\,\hat{Q}_1((p-1)\rho+w_0\lambda)$ is divisible by $\chi((p-1)\rho)$ so we must have  $\psi=s(\lambda)$. So we get $Z=0$ and ${\rm ch}\,\hat{Q}_1((p-1)\rho+w_0\lambda)= {\rm ch}\,({\rm St} \otimes L(\lambda))$. This proves that $${\rm St} \otimes L(\lambda) \simeq \hat{Q}_1((p-1)\rho+w_0\lambda).$$

\noindent Now by \cite[4.2, Satz]{Jantzen 1979}, $\hat{Q}_1((p-1)\rho+w_0\lambda)$ is indecomposable as $G_1$ module, so ${\rm St} \otimes L(\lambda)$ is indecomposable as $G_1$ module. Hence ${\rm St} \otimes L(\lambda)|_{G_{1}}$ is indecomposable. This proves part (b).

\bigskip

\noindent Since ${\rm St} \otimes L(\lambda)|_{G_{1}}$ is indecomposable by \cite[2.1]{Donkin 1993} we get $$ {\rm St} \otimes L(\lambda) \otimes T(\mu)^{[r]} \simeq T((p-1)\rho+\lambda+p^{r}\mu).$$ This gives us result in part (c).

\bigskip

\noindent {\bf Proposition 2.} Suppose $\lambda$ is $r-$minuscule and $\mu \in X^{+}(T)$ then $${\rm St}_{r} \otimes L(\lambda) \otimes T(\mu)^{[r]} \simeq T((p^{r}-1)\rho+\lambda+p^{r}\mu).$$

\bigskip

\noindent {\bf Proof.} Using Steinberg's tensor product theorem we get $$ {\rm St}_{r} \otimes L(\lambda) \simeq \bigotimes_{j=0}^{r-1}\,({\rm St} \otimes L(\lambda^{j}))^{[j]}$$ where $\lambda$ is $r-$minuscule. By above remark we have $$ {\rm St}_{r} \otimes L(\lambda) \simeq \bigotimes_{j=0}^{r-1}\,(T((p-1)\rho+\lambda^{j}))^{[j]}.$$ Apply lemma 1(c) inductively to get $${\rm St}_{r} \otimes L(\lambda) \simeq T((p^{r}-1)\rho+\lambda).$$ Now tensor both sides by $T(\mu)^{[r]}$ and apply lemma 1(c) again to get the result.

\bigskip

\noindent {\bf Corollary.} Let $\lambda$ is $r-$minuscule and $\mu \in X^{+}(T)$ then:

{\bf (a)} $T((p^{r}-1)\rho+p^{r}\mu) \otimes L(\lambda) \simeq T((p^{r}-1)\rho+\lambda+p^{r}\mu).$

{\bf (b)} if $T(\mu)$ is simple then ${\rm St}_{r} \otimes L(p^{r}\mu+\lambda) \simeq T((p^{r}-1)\rho+p^{r}\mu+\lambda)$.

\bigskip

\noindent {\bf Proof.} By lemma 1(a) we get ${\rm St}_{r} \otimes T(\mu)^{[r]} \simeq T((p^{r}-1)\rho+p^{r}\mu)$. Tensor this with $L(\lambda)$ to get the result in part (a).

\bigskip

\noindent If $T(\mu)$ is simple then $L(\mu)\simeq T(\mu)$. So  $L(\lambda) \otimes T(\mu)^{[r]} \simeq L(\lambda) \otimes L(\mu)^{[r]}$. Using Steinberg's tensor product theorem we get $L(\lambda)\otimes L(\mu)^{[r]} \simeq L(\lambda+p^{r}\mu)$. Tensor this with $r-$th Steinberg module to get the result in part (b).

\bigskip

\noindent In case $\lambda$ is $p-$minuscule it is of interest to determine the decomposition ${\rm St} \otimes {\rm L}(\lambda)$, ${\rm St} \otimes \Delta(\lambda)$ and ${\rm St} \otimes \nabla(\lambda)$ as a direct sum of indecomposable modules. In what follows we will show that these are all tilting modules and the direct sum decomposition is determined by the characters of $\nabla(\lambda)$ and ${\rm L}(\lambda)$. We will also show that if $\lambda$ is $(p,r)-$minuscule then ${\rm St_{r}} \otimes {\rm L}(\lambda)$ is tilting. We will also give decomposition of $St_{r} \otimes {\rm L}(\lambda)$ into indecomposable tilting modules.

\bigskip  

\noindent {\bf Lemma 2.} Suppose $\lambda$ is $p-$minuscule. Then every weight $\mu$ of $V(\lambda)$ satisfies $p\rho+\mu \in X^{+}(T)$, where $V(\lambda)=\Delta(\lambda)$ or $\nabla(\lambda)$.

\bigskip

\noindent {\bf Proof.} If $\tau$ is a dominant weight of $V(\lambda)$ then $\tau$ is also $p-$minuscule because $\lambda$ is the dominant weight so $\tau \leq \lambda$ and we can write $\lambda=\tau+\theta$ where $\theta$ is a sum of positive roots. Also $p \geq \langle \lambda,\beta_{0}^{v} \rangle = \langle \tau,\beta_{0}^{v} \rangle+\langle \theta ,\beta_{0}^{v} \rangle \geq \langle \tau,\beta_{0}^{v} \rangle$.

\bigskip

\noindent Let $\mu$ be a weight of $V(\lambda)$ then $w\mu=\tau$ for some $p-$minuscule $\tau \in X^{+}(T)$ and $w \in W$. Let $\alpha$ be a simple root then $\langle p\rho+\mu,\alpha^{v} \rangle = p+\langle w^{-1}\tau,\alpha^{v} \rangle=p+\langle \tau,(w\alpha)^{v} \rangle$. So we need to show that $p+\langle \tau,\gamma^{v} \rangle \geq 0$ for all roots $\gamma$.

\bigskip

\noindent Now  $p+\langle \tau,\gamma^{v} \rangle \geq 0$ for all roots $\gamma$ $\Longleftrightarrow  p+\langle \tau,(w_{0}\gamma)^{v} \rangle \geq 0$ for all roots $\gamma$. And this is true $\Longleftrightarrow p+\langle w_{0}\tau,\gamma^{v} \rangle \geq 0 \Longleftrightarrow p-\langle -w_{0}\tau,\gamma^{v} \rangle \geq 0 \Longleftrightarrow p-\langle \tau,\gamma^{v} \rangle \geq 0$. From the last inequality we get $\langle \tau,\gamma^{v} \rangle \leq p$ and since $\langle \tau,\gamma^{v} \rangle \leq \langle \tau,\beta_{0}^{v} \rangle \leq p$ we have the required result.

\bigskip

\noindent Recall that if $ 0 = M_{0} \leq M_{1} \leq ... \leq M_{t} = M $ is a chain of $B-$modules and $R {\rm ind} _{B}^{G} M_{i} / M_{i-1} = 0, 1\leq i \leq t$ then for ${\rm ind}_{B}^{G} M$ we have a sequence $0= {\rm ind}_{B}^{G} M_{0} \leq {\rm ind}_{B}^{G} M_{1} \leq ... \leq {\rm ind}_{B}^{G} M_{t}={\rm ind}_{B}^{G} M$ with ${\rm ind}_{B}^{G} M_{i}/{\rm ind}_{B}^{G} M_{i-1} \simeq {\rm ind}_{B}^{G} M_{i}/M_{i-1}$. This follows by induction on $t$. Recall also that $R {\rm ind} _{B}^{G} \mu =0$ if $\langle \mu , \alpha^{v} \rangle \geq -1$ for all simple roots $\alpha$. This follows by Kempf's vanishing theorem and \cite[II, proposition 5.4(a)]{Jantzen}. 

\bigskip

\noindent {\bf Proposition 3.} Assume $\lambda$ is $p-$minuscule and let $V$ be a finite dimensional $G-$module such that $\mu \leq \lambda$ for all weights $\mu$ of $V$. Then ${\rm St} \otimes V$ is a tilting module. 

\bigskip

\noindent {\bf Proof.} We will show that ${\rm St} \otimes V$ has a $\nabla-$ filtration. Let $\mu$ be a weight of $V$, then $\mu$ is a weight of some composition factor $L(\nu)$ of $V$. Now $\nu \leq \lambda$, so $\langle \nu,\beta_{0}^{v} \rangle \leq \langle \lambda,\beta_{0}^{v} \rangle \leq p$, therefore $\nu$ is $p-$minuscule. Moreover $\mu$ is a weight of $L(\nu)$ and hence of $\nabla(\nu)$ and so by lemma 2 we have $p\rho+\mu \in X^{+}(T)$. 

\bigskip

\noindent Now choose a $B-$module filtration of $V$ given by $0= V_{0} \leq V_{1} \leq ... \leq V_{t}=V$ with $V_{i}/V_{i-1} \simeq k_{\mu_{i}}$ where $\mu_{i}$ is a weight of $V$. Then ${\rm St} \otimes V= {\rm ind}_{B}^{G}((p-1)\rho \otimes V)$ and $(p-1)\rho \otimes V$ has a filtration $0= (p-1)\rho \otimes V_{0} \leq (p-1)\rho \otimes V_{1} \leq ... \leq (p-1)\rho \otimes V_{t}=(p-1)\rho \otimes V$.

\bigskip

\noindent Also for each section $(p-1)\rho \, \otimes \, V_{i}/V_{i-1}$ we have $R{\rm ind}_{B}^{G}((p-1)\rho \otimes V_{i}/V_{i-1})=R{\rm ind}_{B}^{G}((p-1)\rho \otimes k_{\mu_{i}})=R{\rm ind}_{B}^{G}((p-1)\rho+\mu_{i})=0$ because $\langle (p-1)\rho+\mu_{i},\alpha^{v} \rangle \geq -1$. So ${\rm St} \otimes V$ has filtration in section 

$${\rm ind}_{B}^{G}((p-1)\rho \otimes V_{i}/V_{i-1})=\begin{cases}
\nabla(\mu_{i}), &\mu_{i} \in X^{+}(T)\\
0,& {\rm otherwise.}
\end{cases}$$   

\noindent Therefore ${\rm St} \otimes V$ has a $\nabla-$filtration. Also $\mu^{*} \leq \lambda^{*}$ for all weights $\mu^{*}$ of $V^{*}$ and $\lambda^{*}$ is $p-$minuscule. So ${\rm St} \otimes V^{*}$ has a $\nabla-$filtration. Therefore $({\rm St} \otimes V^{*})^{*}={\rm St} \otimes V$ has a $\Delta-$filtration. Hence ${\rm St} \otimes V$ is tilting. 

\bigskip

\noindent {\bf Corollary.}  Suppose $\lambda$ is $p-$minuscule then ${\rm St} \otimes \Delta(\lambda) \simeq {\rm St} \otimes \nabla(\lambda)$. 

\bigskip

\noindent {\bf Proof.} By proposition 3, ${\rm St} \otimes \Delta(\lambda)$ and ${\rm St} \otimes \nabla(\lambda)$ are tilting modules. Moreover ${\rm St} \otimes \Delta(\lambda)$ and ${\rm St} \otimes \nabla(\lambda)$ have the same character and hence are isomorphic. 

\bigskip

\noindent {\bf Theorem 1.} Let $\lambda$ is $p-$minuscule and $V$ be a finite dimensional $G-$module such that $\mu \leq \lambda$ for all weights $\mu$ of $V$. Then $${\rm St} \otimes V \simeq \bigoplus_{\nu \in X^{+}(T)}\, a_{\nu} \,T((p-1)\rho+\nu)$$

\noindent where ${\rm ch}(V)=\Sigma_{\nu \in X^{+}(T)} \, a_{\nu} s(\nu)$.

\bigskip

\noindent {\bf Proof.} By proposition 3 we have ${\rm St} \otimes V$ is a tilting module. Also by \cite[proposition 5.5]{Donkin 2007} we get ${\rm ch} \, T((p-1)\rho+\nu)= \chi((p-1)\rho)s(\nu)$. Write ${\rm ch}\,(V)=\Sigma_{\nu \in X^{+}(T)} \, a_{\nu} s(\nu)$ then the tilting modules ${\rm St} \otimes V$ and \newline $\oplus_{\nu \in X^{+}(T)} \, a_{\nu} T((p-1)\rho+\nu)$ have the same character and hence are isomorphic.

\bigskip

\noindent {\bf Proposition 4.} Assume $\lambda$ is $(p,r)-$minuscule then ${\rm St_{r}} \otimes L(\lambda)$ is a tilting module.

\bigskip

\noindent {\bf Proof.} Since $\lambda$ is $(p,r)-$minuscule this implies $\lambda \in X_{r}(T)$ and $\lambda=\Sigma_{j=0}^{r-1}\, p^{j} \lambda^{j}$, where each $\lambda^{j}$ is $p-$minuscule. Using Steinberg tensor product theorem we have ${\rm St_{r}} \otimes L(\lambda)= \bigotimes_{j=0}^{r-1}\,({\rm St} \otimes L(\lambda^{j}))^{[j]}$. By proposition 3, ${\rm St}\otimes L(\lambda^{j})$ is tilting for each $\lambda^{j}$. We will use mathematical induction to complete the proof. 

\bigskip

\noindent Write ${\rm St_{r}} \otimes L(\lambda)= {\rm St} \otimes L(\lambda^{0}) \otimes ({\rm St} \otimes L(\lambda^{1})\otimes {\rm St}^{[1]} \otimes L(\lambda^{2})^{[1]} \otimes ... \otimes {\rm St}^{[r-2]} \otimes L(\lambda^{r-1})^{[r-2]})^{[1]}$. Then using inductive hypothesis and theorem 1 we have ${\rm St_{r}} \otimes L(\lambda)= \bigoplus_{\mu} \, a_{\mu} {\rm St} \otimes L(\lambda^{0}) \otimes T(\mu)^{[1]}$. Also by theorem 1, ${\rm St} \otimes L(\lambda^{0}) =\bigoplus_{\nu \in X^{+}(T)}\, b_{\nu} \,T((p-1)\rho+\nu)$. So ${\rm St_{r}} \otimes L(\lambda)=\bigoplus_{\mu , \nu } \, a_{\mu}b_{\nu} \,T((p-1)\rho+\nu) \otimes T(\mu)^{[1]}$. Hence ${\rm St_{r}} \otimes L(\lambda)$ is tilting. 

\bigskip

\noindent {\bf Theorem 2.} Let $\lambda$ is $(p,r)-$minuscule then $${\rm St_{r}} \otimes L(\lambda) \simeq \bigoplus_{\nu \in X^{+}(T)} \, b_{\nu} \,T((p^{r}-1)\rho+\nu)$$

\noindent where ${\rm ch} \, L(\lambda)=\Sigma_{\nu \in X^{+}(T)} \, b_{\nu} s_{r}(\nu)$.

\bigskip

\noindent {\bf Proof.} ${\rm St_{r}} \otimes L(\lambda)$ is tilting by proposition 4. Also by proposition 1 we have ${\rm ch}  \, T((p^{r}-1)\rho+\nu)= \chi((p^{r}-1)\rho)s_{r}(\nu)$. Write ${\rm ch} \, L(\lambda)=\Sigma_{\nu \in X^{+}(T)} \, b_{\nu}   s_{r}(\nu)$ then the tilting modules ${\rm St}_{r} \otimes L(\lambda)$ and \newline $\bigoplus_{\nu \in X^{+}(T)} \, b_{\nu} T((p^{r}-1)\rho+\nu)$ have the same character and hence are isomorphic.

\bigskip

\noindent \textit{Acknowledgements}. I am very grateful to Stephen Donkin for bringing this problem to my attention and for his valuable remarks.

\noindent Department of Mathematics, University of York, Heslington, York, YO10
5DD, United Kingdom.\\ 
E-mail: mfa501@york.ac.uk


\begin{thebibliography}{50}

\bibitem{Bourbaki 1975}
N. Bourbaki. (1975). Groups et algebras de Lie. Chapitres 7 et 8. Hermann.

\bibitem{Donkin 1980}
S. Donkin. (1980). On a question of Verma, J. Lond. Math. Soc. II. Ser. 21:445--455.

\bibitem{Donkin 1993}
S. Donkin. (1993). On tilting modules for algebraic groups, Math. Z. 212:39--60.

\bibitem{Donkin 2007}
S. Donkin.(2007). Tilting modules for algebraic groups and finite dimensional algebras. Handbook of Tilting Theory. pp. 215--257. London Math. Soc. Lecture Notes Series 332 (Cambridge University Press).

\bibitem{Donkin 2007JA}
S. Donkin. (2007). The cohomology of line bundles on the three-dimensional flag variety. J. Algebra. 307:570--613.

\bibitem{Doty 2009}
S.R. Doty. (2009). Factoring Tilting Modules for Algebraic Groups. Journal of Lie Theory. 19(3):531--535.

\bibitem{Jantzen 1979}
J.C. Jantzen. (1979). Uber Darstellungen hoherer Frobenius-Kerne halbein-facher algebraischer Gruppen. Math. Z. 164: 271--292.

\bibitem{Jantzen}
J.C. Jantzen. (2003). Representations of Algebraic Groups, second ed. Math. Surveys Monogr. vol. 107. Amer. Math. Soc.

\bigskip

\end{thebibliography}
\end{document}